\documentclass{article}

\title{\LARGE \textbf{On Relative Length of Long Paths and Cycles in Graphs}}
\author{Zh.G. Nikoghosyan\footnote{G.G. Nicoghossian (up to 1997)}  }

\begin{document}

\maketitle

\begin{abstract}

Let $G$ be a graph on $n$ vertices, $p$ the order of a longest path and $\kappa$ the connectivity of $G$. In 1989, Bauer, Broersma Li and Veldman proved that if $G$ is a 2-connected graph with $d(x)+d(y)+d(z)\ge n+\kappa$ for all triples $x,y,z$ of independent vertices, then $G$ is hamiltonian. In this paper we improve this result by reducing the lower bound  $n+\kappa$ to $p+\kappa$. \\


\end{abstract}

\section{Introduction}

Throughout this article we consider only finite undirected graphs without loops or multiple edges. The set of vertices of a graph $G$ is denoted by $V(G)$ and the set of edges by $E(G)$.  A good reference for any undefined terms is $\cite{[4]}$. For a graph $G$, we use $n$, $\delta$, $\kappa$ and $\alpha$  to denote the order (the number of vertices), the minimum degree, the connectivity and the independence number of $G$, respectively. If $\alpha\ge k$ for some integer $k$, let $\sigma_k$ be the minimum degree sum of an independent set of $k$ vertices; otherwise we let $\sigma_k=+ \infty$.

Each vertex and edge in a graph can be interpreted as simple cycles of orders 1 and 2, respectively. A graph $G$ is hamiltonian if $G$ contains a Hamilton cycle, i.e. a cycle containing every vertex of $G$. A cycle $C$ of a graph $G$ is said to be dominating if $V(G\backslash C)$ is an independent set. The order of a longest path and a longest cycle in $G$ are denoted by $p$ and $c$, respectively. The difference $p-c$ is called relative length denoted by $diff(G)$. A connected graph $G$ is hamiltonian if and only if $diff(G)=0$, that is $c=p$. It is also easy to see that if $diff(G)\le 1$, that is $c\ge p-1$, then any longest cycle in $G$ is a dominating cycle.

The earliest sufficient condition for a graph to be hamiltonian was developed in 1952 due to Dirac \cite{[6]} in terms of order $n$ and minimum degree $\delta$. \\

\noindent\textbf{Theorem A} \cite{[6]}. Every graph with $\delta\ge \frac{n}{2}$ is hamiltonian.  \\

In 1960, Ore \cite{[11]} improved Theorem A by replacing the minimum degree $\delta$ with the arithmetic mean $\frac{1}{2}\sigma_2$ of two smallest degrees among pairwise nonadjacent vertices.\\

\noindent\textbf{Theorem B} \cite{[11]}. Every graph with $\frac{1}{2}\sigma_2\ge \frac{n}{2}$ is hamiltonian.\\

The analog of Theorem A for dominating cycles was established in 1971 by Nash-Williams \cite{[9]}.\\

\noindent\textbf{Theorem C} \cite{[9]}. If  $G$ is a 2-connected graph with $\delta\ge\frac{n+2}{3}$ then each longest cycle in $G$ is a dominating cycle.\\

In 1980, Bondy \cite{[5]} proved the degree sum version of Theorem C. \\

\noindent\textbf{Theorem D} \cite{[5]}. If  $G$ is a 2-connected graph with $\frac{1}{3}\sigma_3\ge \frac{n+2}{3}$ then each longest cycle in $G$ is a dominating cycle.\\

In 1995, Enomoto, Heuvel, Kaneko and Saito \cite{[7]} improved Theorem D by replacing the conclusion "each longest cycle in $G$ is a dominating cycle" with  $c\ge p-1$.\\

\noindent\textbf{Theorem E} \cite{[7]}. If  $G$ is a 2-connected graph with $\frac{1}{3}\sigma_3\ge \frac{n+2}{3}$ then $c\ge p-1$.\\

Using the original proof \cite{[11]}, Theorem B can be essentially improved  by reducing the lower bound  $\frac{n}{2}$ to $\frac{p}{2}$.\\

\noindent\textbf{Theorem 1}. If $G$ is a connected graph with $\frac{1}{2}\sigma_2\ge \frac{p}{2}$ then $c=p=n$.\\

Theorem E can be improved by a similar way based on a result due to Ozeki and Yamashita \cite{[12]}.\\

\noindent\textbf{Theorem 2}. If  $G$ is a 2-connected graph with $\frac{1}{3}\sigma_3\ge \frac{p+2}{3}$ then $c\ge p-1$.\\

The minimum degree versions of Theorems 1 and 2 follow immediately.\\

\noindent\textbf{Corollary 1}. If $G$ is a connected graph with $\delta\ge \frac{p}{2}$ then $c=p=n$.\\

\noindent\textbf{Corollary 2}. If  $G$ is a 2-connected graph with $\delta\ge \frac{p+2}{3}$ then $c\ge p-1$.\\

We propose a conjecture containing Theorems 1 and 2 as special cases when $\lambda =1$ and $\lambda=2$.\\

\noindent\textbf{Conjecture 1}. If $G$ is a $\lambda$-connected graph with
$$
\frac{1}{\lambda+1}\sigma_{\lambda+1}\ge \frac{p+2}{\lambda+1}+\lambda-2
$$
then $c\ge p-\lambda+1$.\\

The long cycles version of Conjecture 1 can be formulated as follows.\\

\noindent\textbf{Conjecture 2}. If $G$ is a $\lambda$-connected $(\lambda\ge 2)$ graph then
$$
c\ge\min\left\{p-\lambda+2,\lambda\left(\frac{1}{\lambda}\sigma_\lambda-\lambda+2\right)\right\}.
$$

Conjecture 2 for $\lambda=2$ was verified  independently by Bondy \cite{[3]}    (1971), Bermond \cite{[2]} (1976) and Linial \cite{[8]} (1976). \\

\noindent\textbf{Theorem F} \cite{[2]}, \cite{[3]}, \cite{[8]}. If $G$ is a 2-connected graph then either $G$ is hamiltonian or $c\ge \sigma_2$.\\

The minimum degree version of Theorem F was proved in 1952 by Dirac \cite{[6]}.\\

\noindent\textbf{Theorem G} \cite{[6]}. If $G$ is a 2-connected graph then either $G$ is hamiltonian or $c\ge 2\delta$.\\

For $\lambda=3$, Conjecture 2 follows immediately from the main result due to Ozeki and Yamashita \cite{[12]}.  \\

\noindent\textbf{Theorem H} \cite{[12]}. If $G$ is a 3-connected graph then either $c\ge \sigma_3-3$ or $c\ge p-1$.\\

In 1981, the bound $n/2$ in Theorem A was reduced to $(n+\kappa)/3$ for 2-connected graphs. \\

\noindent\textbf{Theorem I} \cite{[10]}. If $G$ is a 2-connected graph with $\delta\ge \frac{n+\kappa}{3}$ then $G$ is hamiltonian.  \\

The degree sum version of Theorem I was established in 1989 due to Bauer, Broersma, Li and Veldman \cite{[1]}. \\

\noindent\textbf{Theorem J} \cite{[1]}. If $G$ is a 2-connected graph with $\frac{1}{3}\sigma_3\ge \frac{n+\kappa}{3}$ then $G$ is hamiltonian.  \\

The main result of this paper can be considered as an improvement of  Theorem I by reducing the bound  $(n+\kappa)/3$ to $(p+\kappa)/3$. \\

\noindent\textbf{Theorem 3}. If $G$ is a 2-connected graph with $\frac{1}{3}\sigma_3\ge \frac{p+\kappa}{3}$ then $c=p=n$.  \\

The minimum degree version of Theorem 3 follows immediately.\\

\noindent\textbf{Corollary 3}. If $G$ is a 2-connected graph with $\delta\ge \frac{p+\kappa}{3}$ then $c=p=n$.\\

The following conjecture contains Theorem 3 as a special case when $\lambda=2$.\\

\noindent\textbf{Conjecture 3}. If $G$ is a $\lambda$-connected $(\lambda\ge 2)$ graph with
$$
\frac{1}{\lambda+1}\sigma_{\lambda+1}\ge \frac{p+\kappa+3}{\lambda+1}+\lambda-3
$$
then $c\ge p-\lambda+2$.\\

The long cycle version of Conjecture 3 can be formulated as follows.\\

\noindent\textbf{Conjecture 4}. If $G$ is a $\lambda$-connected $(\lambda\ge 3)$ graph then either
$$
c\ge \lambda\left(\frac{1}{\lambda}\sigma_\lambda-\frac{\kappa}{\lambda}-\lambda+3\right)
$$
or $c\ge p-\lambda+3$.\\

Conjecture 4 for $\lambda=3$ was verified by Yamashita \cite{[13]}.\\

\noindent\textbf{Theorem K} \cite{[13]}. If $G$ is a 3-connected graph then either $c\ge \sigma_3-\kappa$ or $G$ is hamiltonian.\\

The minimum degree version of Theorem K was established by the author \cite{[10]}.\\

\noindent\textbf{Theorem L} \cite{[10]}. If $G$ is a 3-connected graph then either $c\ge 3\delta-\kappa$ or $G$ is hamiltonian.\\

To prove Theorem 2, we need the following result due to Ozeki and Yamashita \cite{[12]}.\\

\noindent\textbf{Theorem M} \cite{[12]}. If $G$ is a 2-connected graph then either $c\ge p-1$ or $c\ge \sigma_3-3$ or $\kappa=2$ and $p\ge \sigma_3-1$.

\section{Proofs}

First we introduce some additional notation.

If $P$ is a path in a graph $G$ then we denote by $\overrightarrow{P}$ the path $P$ with a given orientation, and by $\overleftarrow{P}$ the same path with reverse orientation. If $u,v\in V(P)$ and $u$ precedes $v$ on $\overrightarrow{P}$ then $u\overrightarrow{P}v$ denotes the consecutive vertices of $P$ from $u$ to $v$. The same vertices in reverse order are given by $v\overleftarrow{P}u$. We will consider $u\overrightarrow{P}v$ and $v\overleftarrow{P}u$ both as paths and as vertex sets. If $u\in V(P)$ then $u^+$ denotes the successor of $u$  on $\overrightarrow{P}$ and $u^-$ its predecessor. For $U\subseteq V(P)$, $U^+=\{u^+|u\in U\}$ and $U^-=\{u^-|u\in U\}$. Similar notation is used for cycles.

The proof of Theorem 1 is based on standard arguments originally proposed by Ore \cite{[11]}.\\

\noindent\textbf{Proof of Theorem 1}. Let $G$ be a connected graph with $\sigma_2\ge p$ and let $\overrightarrow{P}=x\overrightarrow{P}y$ be a longest path in $G$ of order $p$. Clearly,  $N(x)\cup N(y)\subseteq V(P)$. \\

\textbf{Case 1}. $xy\in E(G)$.

If $p<n$ then recalling that $G$ is connected, we can construct a path longer than $P$, a contradiction. Otherwise $p=n$, implying that $c=p=n$.\\

\textbf{Case 2}. $xy\not\in E(G)$.

It follows that $x\not\in N(x)\cup N^+(y)$. If  $N(x)\cap N^+(y)=\emptyset$ then
$$
p\ge |N(x)|+|N^+(y)|+|\{x\}|
$$
$$
=|N(x)|+|N(y)|+1=d(x)+d(y)+1\ge \sigma_2+1,
$$
contradicting the hypothesis. Now let $N(x)\cap N^+(y)\not=\emptyset$ and $z\in N(x)\cap N^+(y)$. Then  $xz\overrightarrow{P}yz^-\overleftarrow{P}x$ is a cycle of order $p$ and we can argue as in Case 1.     \qquad \rule{7pt}{6pt}\\

\noindent\textbf{Proof of Theorem 2}. Let $G$ be a 2-connected graph with $\sigma_3\ge p+2$. By Theorem L, either $c\ge p-1$ or $c\ge \sigma_3-3$ or $\kappa=2$, $p\ge \sigma_3-1$. Recalling that $\sigma_3\ge p+2$ (by the hypothesis), we get either $c\ge p-1$ or $p\ge p+1$. Since the latter is impossible, we have $c\ge p-1$.          \qquad \rule{7pt}{6pt}  \\

\noindent\textbf{Proof of Theorem 3}. Let $G$ be a 2-connected graph with $\sigma_3\ge p+\kappa$. Assume first that $\kappa\ge 3$. By Theorem J, we can assume that $c\ge \sigma_3-\kappa$, implying that $c\ge p$. If $c<n$ then clearly $p\ge c+1$ (since $G$ is connected), contradicting $c\ge p$. Hence $c=p=n$, that is $G$ is hamiltonian.

Now assume that $\kappa=2$. Since $\sigma_3\ge p+\kappa=p+2$, by Theorem 2, $c\ge p-1$, implying that each longest cycle in $G$ is a dominating cycle. Let $C$ be a longest cycle in $G$.\\

\textbf{Case 1}. $d(x)=2$ for some $x\in V(G\backslash C)$.

Since $C$ is a dominating cycle, we have $N(x)\subseteq V(C)$. Set $N_C^+(x)=\{y,z\}$. By the maximality of $C$, we have $xy,xz\not\in E(G)$. We have also $yz\not\in E(G)$, since otherwise
$$
y^-xz^-\overleftarrow{C}yz\overrightarrow{C}y^-
$$
is a cycle longer than $C$. Thus, $\{x,y,z\}$ is an independent set of vertices. Further, if either $N(y)\not\subseteq V(C)$ or $N(z)\not\subseteq V(C)$ then we can form a path of order at least $c+2$, contradicting $c\ge p-1$. Hence, $N(y)\cup N(z)\subseteq V(C)$. Put
$$
A=V(y^+\overrightarrow{C}z),  \  \  B=V(z^+\overrightarrow{C}y).
$$

If $w\in N_A(y)\cap N_A^+(z)$ then
$$
y^-xz^-\overleftarrow{C}wy\overrightarrow{C}w^-z\overrightarrow{C}y^-
$$
is a cycle longer than $C$, a contradiction. Hence $N_A(y)\cap N_A^+(z)=\emptyset$. By a symmetric argument, $N_B^+(y)\cap N_B(z)=\emptyset$. Then
$$
c\ge |N_A(y)|+|N_B^+(y)|+|N_A^+(z)|+|N_B(z)|
$$
$$
=|N_C(y)|+|N_C(z)|=d(y)+d(z)
$$
$$
=d(x)+d(y)+d(z)-2\ge \sigma_3-2\ge p.
$$

\textbf{Case 2}. $d(v)\ge 3$ for each $v\in V(G\backslash C)$.

Let $S=\{v_1,v_2\}$ be a cut set of $G$ and let $H_1,H_2,...,H_t$ be the components of  $G\backslash S$. \\

\textbf{Case 2.1}. $V(C)\subseteq V(H_i)\cup S$ for some $i\in \{1,2,...,t\}$.

Assume without loss of generality that $V(C)\subseteq V(H_1)\cup S$. Let $u_1\in V(H_2)$. Since $u_1\not\in V(C)$, we have $d(u_1)\ge3$. Then for each $u_2\in N(u_1)\backslash \{v_1,v_2\}$, we have $u_1u_2\in E(G)$ and $u_1,u_2\not\in V(C)$. This means that $C$ is not a dominating cycle, a contradiction.\\

 \textbf{Case 2.2}. $V(C)\not\subseteq V(H_i)\cup S$ \ $(i=1,2,...,t)$.

 It follows that $V(C)\cap V(H_i)\not=\emptyset$ and $V(C)\cap V(H_j)\not=\emptyset$ for some distinct $i,j \in \{1,2,...,t\}$, say $i=1$ and $j=2$. Recalling also that $|S|=2$, we conclude that $V(C)\subseteq V(H_1)\cup V(H_2)\cup S$ and $v_1,v_2\in V(C)$. If $t\ge3$ then we can argue as in Case 2.1. Hence $t=2$. Clearly, $C$ consists of two paths $P_1$ and $P_2$ with common end vertices $v_1,v_2$ and
 $$
 V(P_i)\subseteq V(H_i)\cup S \  \  (i=1,2).
 $$
 In other words, $\overrightarrow{C}=v_1\overrightarrow{P_1}v_2\overrightarrow{P_2}v_1$.
 Further, if $V(C)= V(H_1)\cup V(H_2)\cup S$ then $c=p=n$, and we are done.
 Otherwise we can choose $x\in V(G\backslash C)$.
 Since $v_1,v_2\in V(C)$, we have $x\in V(H_i)$ for some $i\in \{1,2\}$,
 say $x\in V(H_1)$. We have $N(x)\subseteq V(C)$,
  since $C$ is a dominating cycle. Choose $y\in N^+(x)$ such that
  $|v_1\overrightarrow{P_1}y|$ is as small as possible. If $w\in N(x)\cap N^-(y)$ then
 $$
 v_1\overrightarrow{C}y^-xw\overleftarrow{C}yw^+\overrightarrow{C}v_1
 $$
is a cycle longer than $C$, a contradiction. Hence, $N(x)\cap N^-(y)=\emptyset$, implying that
$$
|P_1|\ge |N(x)|+|N^-(y)|-|\{v_1^-\}|\ge d(x)+d(y)-1.
$$

\textbf{Case 2.2.1}. $V(P_2)=V(H_2)\cup S$.

Clearly, $|P_2|\ge |N(z)|+|\{z\}|\ge d(z)+1$ for each $z\in V(H_2)$ and $\{x,y,z\}$ is an independent set of vertices. Then
$$
c\ge |P_1|+|P_2|-|\{v_1,v_2\}|
$$
$$
\ge (d(x)+d(y)-1)+(d(z)+1)-2\ge \sigma_3-2\ge p.
$$

\textbf{Case 2.2.2}. $V(P_2)\not=V(H_2)\cup S$.

Let $z\in V(H_2)\backslash V(P_2)$. Since $C$ is a dominating cycle, we have $N(z)\subseteq V(C)$. Then, since $C$ is extreme, $|P_2|\ge |N(z)|+|N^+(z)|-1\ge 2d(z)-1$. Observing also that $\{x,y,z\}$ is an independent set of vertices, we get
$$
c\ge |P_1|+|P_2|-2\ge (d(x)+d(y)-1)+(2d(z)-1)-2
$$
$$
\ge (\sigma_3-2)+d(z)-2\ge \sigma_3-2\ge p.  \  \  \  \  \ \qquad \rule{7pt}{6pt}
$$

\noindent Institute for Informatics and Automation Problems\\ National Academy of Sciences\\
P. Sevak 1, Yerevan 0014, Armenia\\
E-mail: zhora@ipia.sci.am


\begin{thebibliography}{20}

\bibitem{[1]} D. Bauer, H.J. Broersma, R. Li and H.J. Veldman, A generalization of a result of Haggkvist and Nicoghossian, J. Combin. Theory B47 (1989) 237-243.

\bibitem{[2]} J.C. Bermond, On Hamiltonian walks, Congr Numer 15 (1976) 41-51.

\bibitem{[3]} J.A. Bondy, Large cycles in graphs, Discrete Math. 1 (1971) 121-131.

\bibitem{[4]} J.A. Bondy and U.S.R. Murty, Graph Theory with Applications, Macmillan, London and Elsevier, New York (1976).

\bibitem{[5]} J.A. Bondy, Longest paths and cycles in graphs of high degree, Research Report CORR 80-16. University of Waterloo, Waterloo, Ontario, 1980.

\bibitem{[6]} G.A. Dirac, Some theorems on abstract graphs, Proc. London, Math. Soc. 2 (1952) 69-81.

\bibitem{[7]} H. Enomoto, J. van den Heuvel, A, Kaneko and A. Saito, Relative length of long paths and cycles in graphs with large degree sums, J. Graph Theory 20 (1995) 213-225.

\bibitem{[8]} N. Linial, A lower bound on the circumference of a graph, Discrete Math. 15 (1976) 297-300.

\bibitem{[9]} C.St.J.A. Nash-Williams, Edge-disjoint hamiltonian cycles in graphs with vertices of large valency, in: L. Mirsky, ed., "Studies in Pure Mathematics", pp. 157-183, Academic Press, San Diego/London (1971).

\bibitem{[10]} Zh.G. Nikoghosyan, On maximal cycle of a graph, DAN Arm. SSR, v. LXXII, 2 (1981) 82-87.

\bibitem{[11]} O. Ore, A note on hamiltonian circuits, Am. Math. Month. 67 (1960) 55.

\bibitem{[12]} K. Ozeki and T. Yamashita, Length of longest cycles in a graph whose relative length is at least two, Graphs and combinatorics, 28 (2012) 859-868.

\bibitem{[13]} T. Yamashita, A degree sum condition for longest cycles in 3-connected graphs, J. Graph Theory 54 (2007) 277-283.

\end{thebibliography}
\end{document}